\documentclass{article}
\usepackage{latexsym,amsfonts,amsmath,amsthm,makeidx}
\usepackage{makeidx}
\makeindex
\newtheorem{definition}{Definition}[section]
\newtheorem{theorem}{Theorem}[section]
\newtheorem{lemma}{Lemma}[section]
\newtheorem{corollary}{Corollary}[section]
\newtheorem{proposition}{Proposition}[section]
\newtheorem{remark}{Remark}[section]
\newcommand{\RN}{\mathbb R^N}
\newcommand{\om}{\omega}

\newcommand{\iy}{\infty}

\newcommand{\s}{\section}
\newcommand{\dd}{\delta}
\newcommand{\g}{\gamma}
\newcommand{\G}{\Gamma}

\newcommand{\la}{\lambda}

\newcommand{\R}{\mathbb R}
\newcommand{\al}{\alpha}

\newcommand{\bb}{\beta}

\newcommand{\m}{\mathcal{M}}

\newcommand{\e}{\varepsilon}
\newcommand{\vp}{\varphi}

\newcommand{\bt}{\begin{theorem}}
\newcommand{\et}{\end{theorem}}
\newcommand{\bl}{\begin{lemma}}
\newcommand{\el}{\end{lemma}}
\newcommand{\bd}{\begin{definition}}
\newcommand{\ed}{\end{definition}}
\newcommand{\bc}{\begin{corollary}}
\newcommand{\ec}{\end{corollary}}
\newcommand{\bp}{\begin{proof}}
\newcommand{\ep}{\end{proof}}
\newcommand{\bx}{\begin{example}}
\newcommand{\ex}{\end{example}}
\newcommand{\bi}{\begin{exercise}}
\newcommand{\ei}{\end{exercise}}
\newcommand{\bo}{\begin{proposition}}
\newcommand{\eo}{\end{proposition}}
\newcommand{\br}{\begin{remark}}
\newcommand{\er}{\end{remark}}
\newcommand{\be}{\begin{equation}}
\newcommand{\ee}{\end{equation}}
\newcommand{\ba}{\begin{align}}
\newcommand{\ea}{\end{align}}
\newcommand{\bn}{\begin{enumerate}}
\newcommand{\en}{\end{enumerate}}
\newcommand{\bg}{\begin{align*}}
\newcommand{\bcs}{\begin{cases}}
\newcommand{\ecs}{\end{cases}}

\newcommand{\intR}[1]{\int_{\R^N}#1\,  dx}

\newcommand{\intRN}[1]{\int_{\R^3}#1\,  dx}

\newcommand{\sg}{\sigma}
\newcommand{\bean}{\begin{eqnarray*}}
\newcommand{\eean}{\end{eqnarray*}}

\numberwithin{equation}{section}

\begin{document}

\title{\bf Standing waves for coupled nonlinear
  Schr\"{o}dinger equations with decaying potentials\thanks{Supported by NSFC (11025106, 11371212, 11271386) and the Both-Side Tsinghua Fund.  E-mail address:
chenzhijie1987@sina.com(Chen);\quad
\quad wzou@math.tsinghua.edu.cn
(Zou)}}
\date{}
\author{{\bf Zhijie Chen$^{1}$,
  Wenming Zou$^2$}\\
\footnotesize {\it  $^{1, 2}$Department of Mathematical Sciences, Tsinghua
University, Beijing 100084, China}\\
\footnotesize {\it $^{1}$Current address: Center for Advanced Study in Theoretical Science,}\\
\footnotesize{\it National Taiwan
University, Taipei 106, Taiwan}}

\maketitle
\begin{center}
\begin{minipage}{120mm}
\begin{center}{\bf Abstract}\end{center}

We study the following singularly perturbed problem for a coupled nonlinear
Schr\"{o}dinger system which arises in Bose-Einstein condensate: $-\e^2\Delta u +a(x) u =
\mu_1 u^3+\beta uv^2$ and $-\e^2\Delta v +b(x) v =
\mu_1 v^3+\beta u^2v$ in $\R^3$ with $u, v> 0$ and
$u(x), v(x)\to 0$ as $|x|\to \iy$. Here, $a, b$ are nonnegative continuous potentials, and
$\mu_1,\mu_2>0$. We consider the case where the coupling constant $\beta>0$ is relatively large. Then for sufficiently small $\e>0$, we obtain
positive solutions of this system which concentrate around local minima of the potentials as $\e\to 0$.
The novelty is that the potentials $a$ and $b$
may vanish at someplace and decay to $0$ at infinity.
\end{minipage}
\end{center}

\section{Introduction}

In this paper we consider standing wave solutions of time-dependent coupled nonlinear
Schr\"{o}dinger equations:
\be\label{eq1}
\begin{cases}
-i\hbar\frac{\partial}{\partial t}\Phi_1-\frac{\hbar^2}{2}\Delta \Phi_1+ a(x)\Phi_1=
\mu_1 |\Phi_1|^2 \Phi_1+\beta |\Phi_2|^2\Phi_1,\,\, x\in \RN, \,\,t>0,\\
-i\hbar\frac{\partial}{\partial t}\Phi_2-\frac{\hbar^2}{2}\Delta \Phi_2+ b(x)\Phi_2=
\mu_2|\Phi_2|^2 \Phi_2+\beta |\Phi_1|^2\Phi_2,\,\, x\in \RN, \,\,t>0,\\
\Phi_j=\Phi_j(x,t)\in\mathbb{C},\quad j=1,2,\\
\Phi_j(x,t)\to 0,  \quad\hbox{as}\,\, |x|\to+\iy, \,\,t>0, \,\,j=1,2,
\end{cases}\ee
where $N\le 3$, $i$ is the imaginary unit, $\hbar$ is the Plank constant,
$\mu_1,\mu_2 >0$ and $\beta\neq 0$ is a coupling constant.
The system (\ref{eq1}) appears
in  the Hartree-Fock theory for a double condensate, i.e., a binary mixture of Bose-Einstein
condensates in two different hyperfine states $|1\rangle$ and $|2\rangle$ (cf. \cite{EGBB}).
Physically, $\Phi_j$ are the corresponding condensate amplitudes, $\mu_j$ and $\beta$ are the intraspecies and
interspecies scattering lengths. The sign of $\beta$ determines whether the interactions of
states $|1\rangle$ and $|2\rangle$
are repulsive or attractive, i.e., the interaction is attractive if $\beta>0$, and the interaction is repulsive
if $\beta<0$, where the two states are in strong competition.

To obtain standing waves of the system (\ref{eq1}), we set $\Phi_1(x, t)=e^{-i E t/\hbar}u(x)$
and $\Phi_2(x, t)=e^{-i E t/\hbar}v(x)$. Then the system (\ref{eq1})
is reduced to the following elliptic system

\be\label{eq2}
\begin{cases}-\frac{\hbar^2}{2}\Delta u +(a(x)-E) u =
\mu_1 u^3+\beta uv^2, \quad x\in \R^N,\\
-\frac{\hbar^2}{2}\Delta v +(b(x)-E) v =\mu_2 v^3+\beta vu^2,  \quad   x\in
\R^N,\\
u(x), v(x)\to 0 \,\,\hbox{as $|x|\to \iy$}.\end{cases}\ee

In this paper we are concerned with positive
solutions for small $\hbar>0$. For sufficiently small
$\hbar>0$, the standing waves are referred to as semiclassical
states. Replacing $a(x)-E, b(x)-E$ by $a(x), b(x)$ for convenience,
we turn to consider the following system
\be\label{eq3}
\begin{cases}-\e^2\Delta u +a(x) u =
\mu_1 u^3+\beta uv^2, \quad x\in \R^N,\\
-\e^2\Delta v +b(x) v =\mu_2 v^3+\beta vu^2,  \quad   x\in
\R^N,\\
u> 0, v> 0 \,\,\hbox{in $\R^N$},\\
u(x), v(x)\to 0 \,\,\hbox{as $|x|\to \iy$},\end{cases}\ee
where $a, b$ are nonnegative continuous functions.

One of the difficulties in the study of (\ref{eq3}) is that it has semi-trivial solutions of type $(u, 0)$ or $(0, v)$. We
call solutions $(u, v)$ with $u\not\equiv 0$ and $v\not\equiv 0$ by nontrivial vector solutions (cf. \cite{IT}). A solution $(u, v)$
with $u>0$ and $v>0$ is called a positive vector solution.

System (\ref{eq3}) has been intensively studied in the past ten years, see
 Lin and Wei \cite{LW3}, Pomponio \cite{P}, Montefusco, Pellacci and Squassina \cite{MPS} and Ikoma and Tanaka \cite{IT} and references therein.
In \cite{LW3}, Lin and Wei studied (\ref{eq3}) by analyzing least energy nontrivial vector solutions.
When $\beta>0$,
they showed the existence of a least energy nontrivial vector solution for small $\e>0$
under suitable conditions on the behavior of $a(x), b(x)$ as $|x|\to +\iy$. In \cite{MPS}, Montefusco, Pellacci and Squassina studied
the case $\beta>0$. They assume that $a, b$ both have positive infimums and there exists $z\in\RN, r>0$ satisfying
$$\min_{|x-z|<r}a(x)<\min_{|x-z|=r}a(x),\quad \min_{|x-z|<r}b(x)<\min_{|x-z|=r}b(x).$$
Then they showed for small $\e>0$ that (\ref{eq3}) has a non-zero solution $(u_\e, v_\e)$ such that $u_\e+v_\e$ has
exactly one global maximum point in $\{x: |x-z|<r\}$. However, when $\bb>0$ is small, one component of $(u_\e, v_\e)$ converges to
$0$ (see Theorem 2.1 (ii) in \cite{MPS}). In \cite{IT}, Ikoma and Tanaka also considered the case $\beta>0$. When $\bb>0$ is relatively small,
they constructed a family of solutions of (\ref{eq3}) which concentrates to a positive vector solution. We also refer to \cite{LW3, P}
for the study of (\ref{eq3}) when $\bb<0$.

\vskip0.1in

{\it Note that in all works \cite{IT, LW3, MPS, P} mentioned above, they all assumed that $a$ and $b$ are positive bounded away from $0$}. In this paper, we
consider the case where $a, b$ may vanish at someplace and decay to $0$ at infinity. In the sequel we assume that
\begin{itemize}

\item [$({\bf V_1})$]  $\,\,a, b\in C(\RN, \R)$ and
$\inf_{x\in\RN}a(x)\ge 0,\,\, \inf_{x\in\RN}b(x)\ge 0.$

\item [$({\bf V_2})$] $$\liminf_{|x|\to+\iy}a(x)|x|^2\log(|x|)>0,\quad \liminf_{|x|\to+\iy}b(x)|x|^2\log(|x|)>0.$$

\item [$({\bf V_3})$] There exists a bounded open domain $\Lambda$ such that
$$\inf_{x\in\overline{\Lambda}}a(x)=a_0>0,\quad \inf_{x\in\overline{\Lambda}}b(x)=b_0>0.$$

\end{itemize}

To study the concentration phenomena of
solutions for system (\ref{eq3}), the following constant
coefficient problem plays an important role:

\be\label{eq4}
\begin{cases}-\Delta u +a(P) u =
\mu_1 u^3+\beta uv^2, & x\in \RN,\\
-\Delta v +b(P) v =\mu_2 v^3+\beta vu^2,    & x\in
\RN,\\
u > 0, v > 0, & x\in \RN,\\
u(x), v(x)\to 0 \,\,\hbox{as $|x|\to \iy$},\end{cases}\ee
where $P\in \overline{\Lambda}$. Then $a(P), b(P)>0$ are positive constants.
Note that system
(\ref{eq4}) appears as a limit problem after a suitable
rescaling of (\ref{eq3}). The existence and the asymptotic behavior of nontrivial vector solutions of (\ref{eq4}) have
received great interest recently, see \cite{AC1, AC2,BDW, BW, BWW, CZ1, DWW, LW1, LW,  MMP,  NR, NTTV, S, WW1, WW2} for example.
Define $H:=H^1(\RN)\times H^1(\RN)$.
 It is well known that solutions
of (\ref{eq4}) correspond to the critical points of $C^2$ functional $L_P: H\to \R$
given by
{\allowdisplaybreaks
\begin{align}\label{eq5}
L_P(u, v)=&\frac{1}{2}\intR{(|\nabla u|^2+a(P) u^2+|\nabla v|^2+b(P) v^2)}\nonumber\\
&-\frac{1}{4}\int_{\RN}(\mu_1 u^4+2\beta u^2v^2 +\mu_2 v^4)\,dx.
\end{align}
}%
Define the Nehari manifold
{\allowdisplaybreaks
\begin{align}\label{eq11}\mathcal{N}_P:=\Big\{(u, v)\in H\setminus\{ (0,0)\},&\,\,\, \intR{
(|\nabla u|^2+a(P) u^2+|\nabla v|^2+b(P) v^2)}\nonumber\\
-& \intR{(\mu_1 u^4+2\beta u^2v^2+\mu_2 v^4)}=0\Big\}, \end{align}
}%
and a constant
\be\label{eq6}\bb_0:=\max\{\mu_1, \mu_2\}\cdot\max_{x\in \overline{\Lambda}}\left\{\frac{a(x)}{b(x)}, \frac{b(x)}{a(x)}\right\}.\ee
By $(V_1)$ and $(V_3)$, one has that $0<\bb_0<\iy$. With the help of \cite[Theorem 2]{S}, we have the following

\bo\label{prop1} Let $\beta>\beta_0$. Then for any $P\in\overline{\Lambda}$, (\ref{eq4}) has a positive radially symmetric vector
solution $(U_P, V_P)\in H$ which is a mountain-pass type solution and satisfies
\be\label{eq7}m(P):=L_P(U_P, V_P)=\inf_{(u, v)\in \mathcal{N}_P}L_P(u, v).\ee
Moreover, $P\mapsto m(P): \overline{\Lambda}\to \R$ is continuous.\eo

\br\label{remark} We call a nontrivial vector solution $(U, V)$ satisfying (\ref{eq7}) a least energy vector solution. So $(U_P, V_P)$
is a positive least energy vector solution.
By \cite{BS}, we see that for $\bb>0$, any positive solution of (\ref{eq4}) is radially symmetric
with respect to some point $x_0\in\RN$.\er

By Proposition \ref{prop1}, $m(P)$ is well defined and continuous in $\overline{\Lambda}$. Assume that

\begin{itemize}

\item [$({\bf V_4})$] There exists a bounded smooth open domain $O\subset\Lambda$ such that
 $$\quad m_0:=\inf\limits_{P\in O} m(P)<\inf\limits_{P\in \partial O}m(P).$$

\end{itemize}

\br\label{remark1}Assumption $(V_4)$ is an abstract condition,
since we can not write down explicitly the function $m(P)$. Such
a type of abstract assumptions for system (\ref{eq4}) can be seen
in \cite{IT, MPS}. This is also a general condition. In the special case of $a(x)=b(x)+C$, where $C\ge 0$ is
a constant, one can easily show that $(V_4)$ holds if
$\inf_{P\in O} a(P)<\inf_{P\in \partial O}a(P).$
Further comments about assumption $(V_4)$ can be seen in \cite[Remarks 1.4-1.5]{IT}.\er

Define
\be\label{eq8}\mathcal{M}:=\left\{P\in O \,\,:\,\, m(P)=m_0\right\}.\ee
Now we can state our main result.

\bt\label{th1}
Let $N=3$, $\bb>\bb_0$ and assumptions $(V_1)-(V_4)$ hold. Then there exists
$\e_0>0$, such that for any
$\e\in(0, \e_0)$, there exists a positive vector solution
$(\tilde{u}_\e, \tilde{v}_\e)$ of (\ref{eq3}), which satisfies
\begin{itemize}
\item [$({i})$] there exists a maximum point $\tilde{x}_\e$ of $\tilde{u}_{\e}+
\tilde{v}_{\e}$ such that $$\lim\limits_{\e\to
0}\hbox{dist}(\tilde{x}_{\e}, \mathcal{M})=0.$$

\item [$({ii})$] for any such $\tilde{x}_{\e}$, $(w_{1,\e}(x), w_{2,\e}(x))
=(\tilde{u}_{\e}(\e x+\tilde{x}_{\e}),\tilde{v}_{\e}(\e
x+\tilde{x}_{\e}))$ converge (up to a subsequence) to a positive least energy vector
solution $(w_1(x), w_2(x))$ of (\ref{eq4}) with $P=P_0$, where $\tilde{x}_{\e} \to P_0 \in
\mathcal{M}$ as $\e \to 0$.

\item [$({iii})$] For any $\al>0$, there exists $c, C>0$ independent of $\e>0$ such that $$(\tilde{u}_{\e}+
\tilde{v}_{\e})(x)\le C \exp\left({-\frac{c}{\e}\frac{|x-\tilde{x}_{\e}|}{1+|x-\tilde{x}_{\e}|}}\right)(1+|x-\tilde{x}_{\e}|)^{-1}\big|\log(2+|x-\tilde{x}_\e|)\big|^{-\al}.$$

\end{itemize}\et

\br\label{remark2}Since the potentials $a, b$ satisfy $(V_1)-(V_2)$, in our proof of Theorem \ref{th1} we need to use the following
Hardy inequality
\be\label{eq9} \frac{(N-2)^2}{4}\intR{\frac{u^2}{|x|^2}} \le\intR{|\nabla u|^2}\quad \forall u\in C_0^{\iy}(\RN),\ee
which holds for $N\ge 3$. On the other hand, problems (\ref{eq3}) and (\ref{eq4}) become critical or supercritical when $N\ge 4$. Here we only consider the subcritical case, so we assume $N=3$ in Theorem \ref{th1}.\er

\br\label{remark3} Assumption $(V_2)$ implies that neither $a$ nor $b$ have compact supports.
For the scalar case (see (\ref{eq10}) below),
$(V_2)$ was introduced
by Bae and Byeon \cite{BB}.
Using Theorem 2 and Theorem 3 from Bae and Byeon \cite{BB}, it is easily seen that,
if $$\text{both}\quad\limsup_{|x|\to+\iy}a(x)|x|^2\log(|x|)=0\quad\text{and}\quad \limsup_{|x|\to+\iy}b(x)|x|^2\log(|x|)=0,$$
system (\ref{eq3}) has no nontrivial $C^2$ solutions for any $\e>0$. This is the reason that we assume $(V_2)$ in Theorem \ref{th1}.\er

For the scalar equation
\be\label{eq10}-\e^2\Delta u +a(x) u=|u|^{p-1}u,\quad x\in\RN,\ee
where $1<p<\frac{N+2}{N-2}$, there are many works on the existence of solutions which concentrate and develop spike layers,
peaks, around some points in $\RN$ while vanishing elsewhere as $\e\to0$. For the case where $\inf_{x\in\RN} a(x)>0$, we refer to
\cite{BJ,DF, DF1} and references therein. For the case where $\inf_{x\in\RN} a(x)=0$, we refer to \cite{AFM, AMR,BB, BW1, BW2, MVS, YZ}
and references therein.

For such a type of system (\ref{eq3}), as far as we know, there is no result on the case
of potentials vanishing at someplace or decaying to $0$ at
infinity, and Theorem \ref{th1} seems to be the first result on this aspect.
This paper is inspired by \cite{BB}, however the method of their proof cannot work here
because of our general
assumption $(V_4)$. In fact, if we assume
\begin{itemize}
\item [$({\bf V'_4})$] There is a bounded open domain $O\subset\Lambda$ and $x_0\in O$ such
that
$$a(x_0)=\inf\limits_{x\in
O}a(x)<\inf\limits_{x\in\partial O}a(x),\quad
b(x_0)=\inf\limits_{x\in O}b(x)<\inf\limits_{x\in\partial O}b(x),$$
\end{itemize}
instead of $(V_4)$ in Theorem \ref{th1}, it might be possible to
prove Theorem \ref{th1} by following Bae and Byeon'
approach in \cite{BB}. It is easy to check that $(V'_4)$ implies $(V_4)$ but the inverse does not hold,
so $(V_4)$ is a more general assumption. Here we will prove Theorem \ref{th1} by developing
further the methods in \cite{BB, MVS}. The approach in \cite{MVS} was developed from \cite{DF}. Remark that the approach in \cite{MVS} cannot work directly in our paper, since, by their approach, it seems that
one can only get the following decay estimate
$$\tilde{u}_{\e}(x)+\tilde{v}_{\e}(x)\le C \exp\left({-\frac{c}{\e}\frac{|x-\tilde{x}_{\e}|}{1+|x-\tilde{x}_{\e}|}}\right)
(1+|x-\tilde{x}_{\e}|)^{-1},$$
which is not enough for us to show that $(\tilde{u}_{\e}, \tilde{v}_{\e})$ is a solution of (\ref{eq3}) (because, in our following proof, $(\tilde{u}_{\e}, \tilde{v}_{\e})$
is obtained as a solution of a modified problem but not as a solution of the original problem (\ref{eq3})).

The rest of this paper proves Theorem \ref{th1}, and we give some notations here. Throughout this paper,
we denote the norm of $L^p(\R^3)$ by $|u|_p =
(\int_{\R^3}|u|^p\,dx)^{\frac{1}{p}}$, and the norm of $H^1(\R^3)$ by $\|u\|=\sqrt{|\nabla u|^2_2+|u|_2^2}$.
We denote positive constants (possibly different) by
$C, c$, and $B(x,r):=\{y\in \RN : |x-y|<r\}$.

\vskip0.2in

\s{The constant coefficient problem}
\renewcommand{\theequation}{2.\arabic{equation}}

In this section, we study the constant coefficient problem (\ref{eq4}) and prove Proposition \ref{prop1}. We assume $N\le 3$ here.
First we recall a result from \cite{S} about the following problem

\be\label{eq2-1}
\begin{cases}-\Delta u + u =
\mu_1 u^3+\beta uv^2, & x\in \RN,\\
-\Delta v +\la v =\mu_2 v^3+\beta vu^2,    & x\in
\RN,\\
u > 0, v > 0, & x\in \RN,\\
u(x), v(x)\to 0 \,\,\hbox{as $|x|\to \iy$},\end{cases}\ee
where $\mu_1,\mu_2,\la>0$. We denote $L_P(u, v), \mathcal{N}_P$ by $L(u, v), \mathcal{N}$
respectively when $a(P), b(P)$ are replaced by $1, \la$ respectively.
Then we have
\bt\label{th2}(see \cite[Theorem 2(iv) and Subsection 3.4]{S}) Assume that $\la\ge 1$. Then for
$$\beta>\max\left\{\mu_1 \la, \,\,\mu_2 \la^{\frac{N}{2}-1}\right\},$$
problem (\ref{eq2-1}) has a positive least energy vector
solution $(U, V)\in H$ which is a radially symmetric mountain-pass type solution and satisfies
$$L(U, V)=\inf_{(u, v)\in \mathcal{N}}L(u, v).$$
Moreover,
$$L(U, V)<\min\left\{\inf_{(u, 0)\in \mathcal{N}} L(u, 0),\quad \inf_{(0, v)\in \mathcal{N}} L(0, v)\right\}. $$ \et

Theorem \ref{th2} implies the following corollary immediately.

\bc\label{corollary1} Let $\beta>\beta_0$, where $\bb_0$ is defined in (\ref{eq6}).
Then for any $P\in\overline{\Lambda}$, (\ref{eq4}) has a positive least energy vector
solution $(U_P, V_P)\in H$ which is a radially symmetric mountain-pass type solution and satisfies
(\ref{eq7}) and
\be\label{eq2-2}m(P)<\min\left\{\inf_{(u, 0)\in \mathcal{N}_P} L_P(u, 0),\quad \inf_{(0, v)\in \mathcal{N}_P} L_P(0, v)\right\}.\ee \ec

Define
{\allowdisplaybreaks
\begin{align}\label{eq2-11}\mathcal{S}(P):=\{(u, v)\in H: &\,\, L'_P(u, v)=0, \,\,L_P(u, v)=m(p), \nonumber\\
& u>0, v>0,\,\, u, v \,\,\hbox{ are radially symmetric}\},\\
\mathcal{S}:=\{(P, u, v)\in \RN &\times H \,\,:\,\, P\in\overline{\Lambda}, \,\,(u, v)\in \mathcal{S}(P)\}.\nonumber\end{align}
}%

Let $(u, v)\in H$ be any a nonnegative solution of (\ref{eq4}) with $L_P(u, v)=m(P)$. Then (\ref{eq2-2}) implies that
$u\not\equiv 0$ and $v\not\equiv 0$. Therefore we have $u>0$ and $v>0$ by the strong maximum principle. By Remark \ref{remark} there
exists some $x_0\in \RN$ such that $(u(\cdot-x_0), v(\cdot-x_0))\in \mathcal{S}(P)$. We have the following properties.

\bl\label{lemma1}
\begin{itemize}
\item [$({i})$] There exists $C_0, C_1,C_2, C_3>0$ such that for all $(P, u, v)\in\mathcal{S}$, there hold
{\allowdisplaybreaks
\begin{gather}\label{eq2-1-1}
    \|u\|, \|v\|\le C_0,\\
    \label{eq2-1-2} |u|_4, |v|_4\ge C_1,\\
    \label{eq2-1-3} u(x),\, v(x),\, |\nabla u(x)|,\, |\nabla v(x)|\le C_2 e^{-C_3|x|}\quad\forall\,x\in\RN.
\end{gather}
}%

\item [$({ii})$] $\mathcal{S}$ is compact in $\RN\times H$.

\item [$({iii})$] $m(P): \overline{\Lambda}\to \R$ is continuous.

\end{itemize}\el

\noindent{\bf Proof. } The proof is something standard. From (\ref{eq7}) it is standard to see that
\be\label{eq2-3}m(P)=\inf_{(u, v)\in \mathcal{N}_P}L_P(u, v)=\inf_{(u, v)\in H\backslash \{(0,0)\}}\max_{t>0}L_P(tu, tv).\ee

(i) Let $a_1=\max_{x\in\overline{\Lambda}}a(x)$. Then it is well known that
$$-\Delta u +a_1 u= \mu_1 u^3,\quad u\in H^1(\RN)$$
has a positive solution $U_0$ which is unique up to a translation. Then
{\allowdisplaybreaks
\begin{align*}
\max_{t>0}L_P(t U_0, 0) &\le \max_{t>0}\left(\frac{1}{2}t^2\intR{(|\nabla U_0|^2+a_1 U_0^2)}-\frac{1}{4}t^4\intR{\mu_1 U_0^4}\right)\\
&=\frac{1}{4}\intR{(|\nabla U_0|^2+a_1 U_0^2)}.
\end{align*}
}%
Combining this with (\ref{eq2-3}) one has that $m(P)$ is uniformly bounded for $P\in \overline{\Lambda}$. Since for any $(P, u, v)\in \mathcal{S}$,
$$4m(P)=\intR{(|\nabla u|^2+a(P) u^2+|\nabla v|^2+b(P) v^2)},$$
we see from $(V_3)$ that (\ref{eq2-1-1}) holds. Recall that for any $(P, u, v)\in \mathcal{S}$, $u, v$
are radially symmetric, so we see from \cite[Lemma A.II]{BL} that
$$u(x),\, v(x)\to 0\quad \hbox{as $|x|\to+\iy$, uniformly for $(P, u, v)\in \mathcal{S}.$}$$
Using a comparison principle, we see that (\ref{eq2-1-3}) holds.
To prove (\ref{eq2-1-2}), we assume by contradiction that there exists a sequence $(P_n, u_n, v_n)\in \mathcal{S}$ such that
\be\label{eq2-4}\lim_{n\to+\iy}|u_n|_4=0.\ee
(The case $|v_n|_4\to 0$ is similar.) Passing to a subsequence, $P_n\to P_0\in\overline{\Lambda}$.
Define
$$H^1_r (\RN):=\{u\in H^1(\RN)\,\,:\,\, u \hbox{ is radially symmetric}\}.$$
Since $L_{P_n}'(u_n, v_n)=0$ and the Sobolev embedding $H^1_r(\RN)\hookrightarrow L^4(\RN)$ is compact, it is standard to show that
$(u_n, v_n)$ converges to some $(u_0, v_0)$ strongly in $H$ (up to a subsequence), $L_{P_0}'(u_0, v_0)=0$ and
\be\label{eq2-5}\lim_{n\to\iy}m(P_n)=\lim_{n\to\iy} L_{P_n} (u_n, v_n)=L_{P_0}(u_0, v_0).\ee
By (\ref{eq2-4}), we get from $L'_{P_n}(u_n, v_n)(u_n, 0)=0$ that
$$\intR{(|\nabla u_n|^2+a(P_n) u_n^2)}=\intR{(\mu_1 u_n^4+\beta u_n^2 v_n^2)}\to 0,\,\,\hbox{as $n\to\iy$,}$$
which implies that $u_0=0$.
We denote $L_P(u, v)$ by $L_0(u, v)$
when $a(P), b(P)$ are replaced by $a_0, b_0$. Then by a standard mountain-pass argument, there exists some $\rho,\al>0$ such that
$\inf\limits_{\|u\|+\|v\|=\rho}L_0(u, v)=\al>0.$ By (\ref{eq2-3}) and $(V_3)$ this means that
$$m(P)\ge\inf_{(u, v)\in H\backslash \{(0,0)\}}\max_{t>0}L_0(tu, tv)\ge\al>0,\quad\forall \,P\in \overline{\Lambda}.$$
Therefore, $v_0\not\equiv 0$. By (\ref{eq2-2}) we have $m(P_0)<L_{P_0}(0, v_0)$. On the other hand, let $(U, V)\in \mathcal{S}(P_0)$, then
$L_{P_0}(U, V)=m(P_0)$. Note that
\begin{align*}
m(P_n)&\le\max_{t>0}L_{P_n}(tU, tV)=\frac{(\intR{(|\nabla U|^2+a(P_n)U^2+|\nabla V|^2+b(P_n)V^2)})^2}
{4\intR{(\mu_1 U^4+2\bb U^2 V^2+\mu_2 V^4)}}\\
&\to \frac{1}{4}\intR{(|\nabla U|^2+a(P_0)U^2+|\nabla V|^2+b(P_0)V^2)}=m(P_0)
\end{align*}
as $n\to\iy$, that is,
\be\label{eq2-6}\lim_{n\to\iy}m(P_n)\le m(P_0)<L_{P_0}(0, v_0)=L_{P_0}(u_0, v_0),\ee
a contradiction with (\ref{eq2-5}). Hence, (\ref{eq2-1-2}) holds.

(ii) For any sequence $(P_n, u_n, v_n)\in \mathcal{S}$, similarly as in the proof of (i), up to a subsequence, we may
assume that $P_n\to P_0$, $(u_n, v_n)\to (u_0, v_0)$ strongly in $H$ and $(u_0, v_0)$ is a
nontrivial vector solution of (\ref{eq3}) with $P=P_0$.
By (\ref{eq7}) we have $L_{P_0}(u_0, v_0)\ge m(P_0)$. Meanwhile, (\ref{eq2-5}) and (\ref{eq2-6})
imply $L_{P_0}(u_0, v_0)\le m(P_0)$. That is,
$L_{P_0}(u_0, v_0)= m(P_0)=\lim_{n\to\iy}m(P_n)$. Since $u_n, v_n>0$ are radially symmetric, we also have that $u_0, v_0>0$ are
radially symmetric. Hence, $(P_0, u_0, v_0)\in\mathcal{S}$.

(iii) follows from the proof of (ii). This completes the proof.\hfill$\square$\\

Proposition \ref{prop1} follows directly from Corollary \ref{corollary1} and Lemma \ref{lemma1}. \hfill$\square$

\vskip0.2in

\s{Proof of Theorem \ref{th1}}
\renewcommand{\theequation}{3.\arabic{equation}}

In this section we assume that $N=3$, $\bb>\bb_0$ and assumptions $(V_1)-(V_4)$ hold. Define $a_\e (x)= a(\e x), \,b_\e (x)= b(\e x)$.
To study (\ref{eq3}), it suffices to consider the following system

\be\label{eq3-1}
\begin{cases}-\Delta u +a_\e u =
\mu_1 u^3+\beta uv^2, & x\in \R^3,\\
-\Delta v +b_\e v =\mu_2 v^3+\beta vu^2,    & x\in
\R^3,\\
u > 0, v > 0, & x\in \R^3,\\
u(x), v(x)\to 0 \,\,\hbox{as $|x|\to \iy$},\end{cases}\ee

Let $H^1_{a, \e}$ (resp. $H^1_{b, \e}$) be the completion of
$C_0^\iy (\R^3)$ with respect to the norm
$$\|u\|_{a, \e}=\left(\intRN{|\nabla u|^2+a_\e u^2}\right)^{\frac{1}{2}}\,\,
\left(\hbox{resp.}\,\,\|u\|_{b, \e}=\left(\intRN{|\nabla u|^2+b_\e
u^2}\right)^{\frac{1}{2}}\right).$$ Define $H_\e := H^1_{a,
\e}\times H^1_{b, \e}$ with a norm
$\|(u, v)\|_\e =\sqrt{\|u\|^2_{a, \e}+\|v\|^2_{b, \e}}.$

From now on, for any set $B\subset\R^3$ and $\e, s>0$, we define
$B_\e:=\{x\in \R^3 : \e x\in B\}$, $B^s:= \{x\in\R^3 : dist(x, B)\le s\}$ and $B_\e^s:=(B^s)_\e$.
Without loss of generality, we may assume that $0\in\m$ and $B(0, \rho_0)\subset O\subset B(0, \rho_1)$ for some $\rho_1> \rho_0>0$.
By $(V_3)-(V_4)$ we can choose $\delta\in (0, \rho_0)$ small such that $dist(\m, \R^3\backslash O)\ge 5\delta$ and
\be\label{eq3-2}\inf_{x\in\overline{O}^{5\delta}}a(x)\ge a_0/2>0,\quad \inf_{x\in\overline{O}^{5\delta}}b(x)\ge b_0/2>0.\ee
For $0<\e<\rho_0$ we define $\gamma_\e : [\rho_0/\e, +\iy)\to (0,+\iy)$ by
\be\label{eq3-3}\g_\e(t):=\frac{\e^2}{t^2\log t},\ee
and$$\chi_{O_\e} (x) := \begin{cases} 1 &
\hbox{if}\quad x\in O_\e, \\ 0 & \hbox{if}\quad x\not\in
O_\e,\end {cases}$$
Denote
$F(s, t):=\frac{1}{4}(\mu_1 s^4 + 2\bb s^2 t^2+ \mu_2 t^4)$, and set
\be\label{eq3-4}
F_\e(x, s, t)=\begin{cases} F(s, t) &
\hbox{if}\quad F(s, t)\le\frac{1}{4}\g_\e^2(|x|), \\ \g_\e(|x|)\sqrt{F(s,t)}-\frac{1}{4}\g_\e^2(|x|)
& \hbox{if}\quad F(s, t)>\frac{1}{4}\g_\e^2(|x|).\end {cases}
\ee
Then we have
\be\label{eq3-5}
\nabla_{(s, t)} F_\e(x, s, t)=\begin{cases} (\mu_1 s^3+\bb st^2, \mu_2 t^3+\bb s^2t) &
\hbox{if}\,\, F(s, t)\le\frac{1}{4}\g_\e^2(|x|), \\ \g_\e(|x|)\frac{(\mu_1 s^3+\bb st^2, \mu_2 t^3+\bb s^2t)}{2\sqrt{F(s,t)}}
& \hbox{if}\,\, F(s, t)>\frac{1}{4}\g_\e^2(|x|).\end {cases}
\ee
This means that $F_\e(x,\cdot)\in C^1(\R^2)$ as a function of $(s, t)$. Define a truncated function
\be\label{eq3-6}G_\e(x, s, t):=\chi_{O_\e}(x) F(s, t)+(1-\chi_{O_\e}(x))F_\e(x, s, t).\ee
By the definition of $\bb_0$ in (\ref{eq6}), one has that $\bb>\max\{\mu_1, \mu_2\}$. Then it is easy to see that
{\allowdisplaybreaks
\begin{gather}\label{eq3-7}G_\e(x, s, t)\le F(s, t),\quad\forall \,x\in \R^3,\\
\label{eq3-8}0\le 4G_\e(x, s, t)= \nabla_{(s, t)} G_\e(x, s, t)(s, t),\quad\forall\, x\in O_\e,\\
\label{eq3-9}2G_\e(x,s, t)\le \nabla_{(s, t)} G_\e(x, s, t)(s, t)\le\sqrt{\bb} \g_\e(|x|)(s^2+t^2),\quad\forall\, x\in\R^3\backslash O_\e.
\end{gather}
}%
Define a functional $J_\e: H_\e\to \R$ by
\be\label{eq3-10}J_\e(u, v):=\frac{1}{2}\|u\|_{a,\e}^2+\frac{1}{2}\|v\|_{b,\e}^2-\intRN{G_\e(x, u^+, v^+)}.\ee
Here and in the following, $u^+(x):=\max\{u(x), 0\}$ and so is $v^+$. By the following Hardy inequality in dimension $N=3$
\be\label{eq3-11} \frac{1}{4}\intRN{\frac{u^2}{|x|^2}} \le\intRN{|\nabla u|^2}\quad \forall u\in C_0^{\iy}(\R^3),\ee
it is standard to show that $J_\e$ is well defined and $J_\e\in C^1(H_\e, \R)$. Furthermore, any critical points of $J_\e$ are weak solutions
of the following system
\be\label{eq3-12}
\begin{cases}-\Delta u +a_\e u =
\partial_u G_\e(x, u^+, v^+), & x\in \R^3,\\
-\Delta v +b_\e v =\partial_v G_\e(x,u^+, v^+),    & x\in
\R^3,\\
u(x), v(x)\to 0 \,\,\hbox{as $|x|\to \iy$}.\end{cases}\ee
For each small $\e>0$, we will find a nontrivial solution of (\ref{eq3-12}) by applying mountain-pass argument to $J_\e$.
Then we shall prove that this solution
is a positive vector solution of (\ref{eq3-1}) for $\e>0$ sufficiently small. This idea was first introduced by del Pino and Felmer \cite{DF}.

\bl\label{lemma2}Let $\e\in(0, \e_1)$ be fixed, where $\e_1$ satisfies $$\frac{\sqrt{\bb}\e_1^2}{\log{(\rho_0/\e_1)}}=1/8.$$
For any $c\in\R$, let $(u_n, v_n)\in H_\e$ be a $(PS)_c$ sequence for $J_\e$, that is,
$$J_\e(u_n, v_n)\to c,\quad J_\e'(u_n, v_n)\to 0.$$
Then, up to a subsequence, $(u_\e, v_\e)$ converge strongly in $H_\e$.\el

\noindent{\bf Proof. } Recall the definition of $\g_\e$ in (\ref{eq3-3}) and $B(0, \rho_0)\subset O$. By Hardy inequality (\ref{eq3-11}), we have
\begin{align}\label{eq3-18}
\frac{1}{4}\int_{\R^3\backslash O_\e}&\sqrt{\bb}\g_\e(|x|)(u_n^2+v_n^2)\,dx\le
\frac{\sqrt{\bb}\e^2}{\log{\rho_0/\e}}\frac{1}{4}\intRN{\frac{u_n^2+v_n^2}{|x|^2}}\nonumber\\
&\le \frac{\sqrt{\bb}\e^2}{\log{\rho_0/\e}}\intRN{|\nabla u_n|^2+|\nabla v_n|^2}\le\frac{1}{8}\|(u_n,v_n)\|_\e^2.
\end{align}
Therefore, we deduce from (\ref{eq3-8}) and (\ref{eq3-9}) that
{\allowdisplaybreaks
\begin{align}\label{eq3-13}
c+& o(\|(u_n,v_n)\|_\e)\ge J_\e(u_n, v_n)-\frac{1}{4}J_\e'(u_n, v_n)(u_n, v_n)\nonumber\\
&=\frac{1}{4}\|(u_n,v_n)\|_\e^2+\intRN{\left(\frac{1}{4}\nabla_{(u, v)}G_\e(x, u_n^+, v_n^+)(u_n, v_n)-G_\e(x, u_n^+, v_n^+)\right)}\nonumber\\
&\ge\frac{1}{4}\|(u_n,v_n)\|_\e^2-\frac{1}{4}\int_{\R^3\backslash O_\e}\sqrt{\bb}\g_\e(|x|)(u_n^2+v_n^2)\,dx\nonumber\\
&\ge\frac{1}{8}\|(u_n,v_n)\|_\e^2,
\end{align}
}%
that is, $\|(u_n, v_n)\|_\e\le C$ for all $n\in\mathbb{N}$. Up to a subsequence, we may assume
that $(u_n, v_n)\rightharpoonup (u, v)$ weakly in $H_\e$ and $u_n\to u, v_n\to v$ strongly in $L^4_{loc} (\R^3)$.
Since there exists $\al_0>0$ depending on $\mu_1, \mu_2,\bb$ only, such that
\be\label{eq3-14}max\{\mu_1 s^2+\bb t^2,\,\, \bb s^2+\mu_2 t^2\}\le \al_0 2\sqrt{F(s, t)},\quad \forall\, s, t\in \R,\ee
from (\ref{eq3-5}) we obtain
\be\label{eq3-15}|\partial_u G_\e(x, u^+, v^+)|\le \al_0\g_\e(|x|)|u|,\,
 |\partial_v G_\e(x, u^+, v^+)|\le \al_0\g_\e(|x|)|v|,\,\forall\,x\in\R^3\backslash O_\e.\ee
Then for any $R\ge \rho_1$, we deduce from (\ref{eq3-3}), (\ref{eq3-11}) and (\ref{eq3-15}) that
{\allowdisplaybreaks
\begin{align*}
\limsup_{n\to\iy}&\|u_n-u\|_{a,\e}^2=\limsup_{n\to\iy}\intRN{(\partial_u G_\e(x, u_n^+, v_n^+)-\partial_u G_\e(x, u^+, v^+))(u_n-u)}\\
\le&\limsup_{n\to\iy}\int_{B(0,R/\e)}(\partial_u G_\e(x, u_n^+, v_n^+)-\partial_u G_\e(x, u^+, v^+))(u_n-u)\,dx\\
&+\limsup_{n\to\iy}\int_{\R^3\backslash B(0,R/\e)}2\al_0\g_\e(|x|)(u_n^2+u^2)\,dx\\
\le & \limsup_{n\to\iy}C \left(\int_{B(0,R/\e)}|u_n-u|^4\,dx\right)^{1/4}+C\frac{8\e^2\al_0}{\log{R/\e}}=C\frac{8\e^2\al_0}{\log{R/\e}}.
\end{align*}
}%
Since $R\ge \rho_1$ is arbitrary, we see that $u_n\to u$ strongly in $H^1_{a,\e}$. Similarly, $v_n\to v$ strongly in $H_{b,\e}^1$. This
completes the proof.\hfill$\square$\\

For any $\e\in(0, \e_1)$ fixed, we define
\be\label{eq3-16}c_\e:=\inf_{\g\in\Phi_\e}\sup_{t\in[0, 1]} J_\e (\g(t)),\ee
where $\Phi_\e=\{\g\in C([0,1], H_\e) \,: \,\g(0)=(0,0), \,\,J_\e(\g(1))<0\}$.

\bl\label{lemma3} For any fixed $\e\in(0, \e_1)$, there exists a nontrivial critical point $(u_\e, v_\e)$ of $J_\e$ such that
$u_\e\ge 0, v_\e\ge 0$ and
$J_\e(u_\e, v_\e)=c_\e>0$. Moreover, at least one of $u_\e>0$ and $v_\e>0$ holds. \el

\noindent{\bf Proof. } Recall that $\bb>\bb_0\ge\max\{\mu_1, \mu_2\}$. By (\ref{eq3-4}) one has that
\be\label{eq3-17}F_\e(x, s, t)\le \g_\e(|x|)\sqrt{F(s, t)}\le\frac{1}{2}\sqrt{\bb}\g_\e(|x|)(s^2+t^2).\ee
Then we deduce from $(V_3)$, (\ref{eq3-10}) and (\ref{eq3-18}) that
{\allowdisplaybreaks
\begin{align*}
J_\e(u, v)&\ge \frac{1}{2}\|(u, v)\|_\e^2-\int_{O_\e}F(u, v)\,dx-\frac{1}{2}\sqrt{\bb}\int_{\R^3\backslash O_\e}\g_\e(|x|)(u^2+v^2)\,dx\\
&\ge  \frac{1}{4}\|(u, v)\|_\e^2-C \int_{O_\e}(u^4+v^4)\,dx\\
&\ge  \frac{1}{4}\|(u, v)\|_\e^2-C\|(u, v)\|_\e^4.
\end{align*}
}%
Hence, there exists $r, \al_1>0$ small such that
$$\inf_{\|(u, v)\|_\e=r}J_\e(u, v)=\al_1>0.$$
This implies that $c_\e\ge\al_1>0$. Choose $\phi\in C_0^{\iy}(O_\e)$ such that $\phi\ge 0$ and $\phi\not\equiv 0$.
Then $G_\e(x, \phi^+, \phi^+)=F(\phi, \phi)$, which implies that $J_\e(t\phi, t\phi)\to -\iy$ as $t\to+\iy$.
That is, $J_\e$ has a mountain-pass structure.
By Lemma \ref{lemma2} and Mountain Pass Theorem (\cite{AR-pass}), there exists a nontrivial critical point $(u_\e, v_\e)$ of $J_\e$
such that $J_\e(u_\e, v_\e)=c_\e>0$. Denote $u_\e^-(x):=\max\{-u_\e(x), 0\}$ and so is $v_\e^-$. Then we see from (\ref{eq3-12}) that
$$\intRN{|\nabla u_\e^-|^2+a_\e |u_\e^-|^2}=0,\quad \intRN{|\nabla u_\e^-|^2+b_\e |u_\e^-|^2}=0,$$
which implies that $u_\e, v_\e\ge 0$. By the strong maximum principle, at least one of $u_\e>0$ and $v_\e>0$ holds.
This completes the proof.\hfill$\square$\\

Define
$$\tilde{c}_\e:=\inf_{(u, v)\in \widetilde{\Phi}_\e}\max_{t>0}J_\e(tu, tv),$$
where $\widetilde{\Phi}_\e:=\{(u, v)\in H_\e\backslash\{(0,0)\} : \int_{O_\e}{((u^+)^2+ (v^+)^2)}\,dx>0\}$. Then we have the following lemma.

\bl\label{lemma4} For any $(u, v)\in \widetilde{\Phi}_\e$, there exists a unique $t_{u, v}>0$ such that
\be\label{eq3-19} J_\e(t_{u, v}u, t_{u, v}v):=\max_{t>0}J_\e(tu, tv).\ee
Moreover, $c_\e=\tilde{c}_\e$.
\el

\noindent{\bf Proof. } Fix any $(u, v)\in \widetilde{\Phi}_\e$. By the definition (\ref{eq3-6}) of $G_\e(x, u, v)$, we have that for any $x\in \R^3$,
$\frac{1}{t}\frac{d}{d t}G_\e(x, t u^+, t v^+)$ is nondecreasing as $t>0$ increases. Moreover, if $x\in O_\e$ and $(u^+)^2(x)+(v^+)^2(x)>0$, then
$\frac{1}{t}\frac{d}{d t}G_\e(x, t u^+, t v^+)$ is strictly increasing as $t>0$ increases. This means that
there exists a unique $t_{u,v}>0$ such that $\frac{d}{dt}J_\e(tu, tv)|_{t_{u, v}}=0$, that is,
$$\|u\|_{a, \e}^2+\|v\|_{b, \e}^2=\intRN{\frac{1}{t_{u, v}}\frac{d}{d t}G_\e(x, t u^+, t v^+)|_{t_{u, v}}}.$$
Since
$$J_\e(tu, tv)\le \frac{1}{2}t^2\|(u, v)\|_\e^2-t^4\int_{O_\e}F(u^+, v^+)\,dx\to-\iy\,\,\,\hbox{as $t\to+\iy$,}$$
we see that (\ref{eq3-19}) holds and
$c_\e\le \max_{t>0}J_\e(tu, tv).$ Therefore, $c_\e\le \tilde{c}_\e$. Meanwhile, since $(u_\e, v_\e)\in\widetilde{\Phi}_\e$,
$c_\e=J_\e(u_\e, v_\e)$ and $t_{u_\e, v_\e}=1$, we have $c_\e\ge \tilde{c}_\e$. This completes the proof.\hfill$\square$

\bl\label{lemma5} Let $(u_\e, v_\e)$ be in Lemma \ref{lemma3}. Then $\limsup\limits_{\e\to0} c_\e\le \inf\limits_{P\in O}m(P)$, and
there exists $C>0$ independent of $\e>0$, such that $\|(u_\e, v_\e)\|_\e^2\le C$.\el

\noindent{\bf Proof. } Fix any $P\in O$ and let $(U_P, V_P)$ be in Corollary \ref{corollary1}.
Take $T>0$ such that $L_P(TU_P, TV_P)\le -1$. Note that there exists $R>0$ such that $B(P, R):=\{x: |x-P|<R\}\subset O$, we take $\phi\in C_0^1(B(0, R), \R)$ with $0\le \phi\le 1$ and $\phi(x)\equiv 1$ for $|x|\le R/2$. Define $\phi_\e(x):=\phi(\e x)$, then $\phi_\e (x-P/\e)\neq 0$ implies $x\in O_\e$. Combining this with (\ref{eq3-6}) and the Dominated Convergence Theorem, one has
{\allowdisplaybreaks
\begin{align*}
&\quad J_\e (t(\phi_\e U_P)(\cdot-P/\e), \,\, t(\phi_\e V_P)(\cdot-P/\e))\\
&=\frac{t^2}{2}\int_{|x|\le R/\e}{(|\nabla(\phi_\e U_P)|^2+a(\e x+P)\phi_\e^2U_P^2)}\\
&\quad+\frac{t^2}{2}\int_{|x|\le R/\e}{(|\nabla(\phi_\e V_P)|^2+b(\e x+P)\phi_\e^2V_P^2)}-\int_{|x|\le R/\e}{F(t\phi_\e U_P, t\phi_\e V_P)}\\
&\to L_P(tU_P, tV_P),\quad\hbox{as $\e\to 0$, uniformly for $t\in [0, T].$}
\end{align*}
}%
So $J_\e (T(\phi_\e U_P)(\cdot-P/\e),  T (\phi_\e V_P)(\cdot-P/\e))<0$ for $\e>0$ sufficiently small.
Since $((\phi_\e U_P)(\cdot-P/\e),(\phi_\e V_P)(\cdot-P/\e))\in \widetilde{\Phi}_\e$, we see from Lemma \ref{lemma4} that
{\allowdisplaybreaks
\begin{align*}
\limsup\limits_{\e\to0} c_\e&\le\limsup\limits_{\e\to0}\max_{t\in[0,T]} J_\e (t(\phi_\e U_P)(\cdot-P/\e),  t (\phi_\e V_P)(\cdot-P/\e))\\
&=\max_{t\in[0,T]}L_P(tU_P, tV_P)=L_P(U_P, V_P)=m(P).
\end{align*}
}%
Since $P\in O$ is arbitrary, we have $\limsup\limits_{\e\to0} c_\e\le \inf\limits_{P\in O}m(P)$. From (\ref{eq3-13}), there exists $C>0$ independent of $\e>0$, such that $\|(u_\e, v_\e)\|_\e^2\le C$. This completes the proof.
\hfill$\square$

\bl\label{lemma6} Let $(u_\e, v_\e)$ be in Lemma \ref{lemma3}. Then there exists $0<\e_2\le\e_1$, such that for any $\e\in (0, \e_2)$,
there holds
$$\|u_\e+v_\e\|_{L^{\iy}(O_\e)}\ge\sqrt{\min\{a_0, b_0\}/\bb}.$$\el

\noindent{\bf Proof. } Without loss of generality, we assume that $a_0\le b_0$. Recalling $\al_0$ in (\ref{eq3-15}), we take $\e_2\le\e_1$ such that
$\frac{\al_0 \e_2^2}{\log \rho_0/\e_2}\le 1/4$. Assume that $\|u_\e+v_\e\|_{L^{\iy}(O_\e)}\le\sqrt{a_0/\bb}$ for some $\e\in (0,\e_2)$, then
$$\mu_1u_\e^3+\bb u_\e v_\e^2+\bb u_\e^2v_\e+\mu_2 v_\e^3\le a_0(u_\e+v_\e)\,\,\,\hbox{in $O_\e$.}$$
Combining this with (\ref{eq3-12}) and (\ref{eq3-15}) we obtain that
\begin{align*}
-\Delta (u_\e+v_\e)+a_\e u_\e+b_\e v_\e\le \chi_{O_\e}a_0(u_0+v_0)+\al_0(1-\chi_{O_\e})\g_\e(|x|)(u_\e+v_\e),
\end{align*}
that is,
$$-\Delta (u_\e+v_\e)\le (1-\chi_{O_\e})\frac{\al_0 \e^2}{\log{\rho_0/\e}}\frac{u_\e+v_\e}{|x|^2},$$
which implies from (\ref{eq3-11}) that
$$\intRN{|\nabla (u_\e+v_\e)|^2}\le \frac{\al_0 \e^2}{\log{\rho_0/\e}}\intRN{\frac{|u_\e+v_\e|^2}{|x|^2}}<\intRN{|\nabla (u_\e+v_\e)|^2}, $$
a contradiction. This completes the proof.
\hfill$\square$

\bl\label{lemma7} Let $(u_\e, v_\e)$ be in Lemma \ref{lemma3}. Let $(\e_n)_{n\ge 1}$ be a sequence with $\e_n\to0$.
Let $k\ge 1$ and for $i\in [1, k]\cap\mathbb{N}$, there is $\{P_n^i\}_{n\ge1}\subset O_{\e_n}$ with $\e_n P_n^i\to P^i\in \overline{O}$ as $n\to\iy$.
If
\begin{gather*}
    \liminf_{n\to\iy}(u_{\e_n}+v_{\e_n})(P_n^i)>0,\,\,\,\forall\,i;\quad
    \lim_{n\to\iy}|P_n^i-P_n^j|=+\iy,\,\,\forall\,\,i\neq j.
\end{gather*}
Then $\liminf\limits_{n\to\iy} c_{\e_n}\ge \sum_{i=1}^k m(P^i)$.\el

\noindent{\bf Proof. } The proof is inspired by \cite{MVS}. For $i\in\{1,\cdots,k\}$, we define
$$(u_n^i, v_n^i):=(u_{\e_n}(\cdot+P_n^i), v_{\e_n}(\cdot+P_n^i)).$$
By Lemma \ref{lemma5}, $(u_{\e_n}, v_{\e_n})$ is uniformly bounded in $H_\e$, so $u_n^i, v_n^i$ are uniformly bounded
in $H^1_{loc}(\R^3)$.  By the system and the elliptic regularity, it is standard to show that $u_n^i, v_n^i$ are uniformly bounded
in $W^{2,q}_{loc}(\R^3)$ for any $q\ge 2$. By the compactness of Sobolev embedding, passing to a subsequence, we may assume that
$u_n^i\to u^i, v_n^i\to v^i$ strongly in $C^1_{loc}(\R^3)$. Moreover, we have $u^i, v^i\ge 0$
 and $u^i(0)+v^i(0)>0$. By Fatou Lemma, for any $R>0$, we have
{\allowdisplaybreaks
\begin{align*}
\int_{B(0, R)}|\nabla u^i|^2+a(P^i)|u^i|^2\,dx
\le &\liminf_{n\to\iy}\int_{B(0, R)}|\nabla u_n^i|^2+a(\e_n x+\e_n P_n^i)|u_n^i|^2\,dx\\
\le&\liminf_{n\to\iy}\int_{\R^3}|\nabla u_{\e_n}|^2+a_{\e_n}|u_{\e_n}|^2\,dx\le C.
\end{align*}
}%
So $u^i\in H^1(\R^3)$. Similarly, $v^i\in H^1(\R^3)$.
Since $O$ is smooth, up to a subsequence, we may assume that $\chi_{O_{\e_n}}(\cdot+P_n^i)$ converges almost everywhere
to $\chi^i$, where $0\le \chi^i\le1$. In fact, $\chi^i$ is either the characteristic function of $\R^3$ or the characteristic function
of a half space. Then it is easy to see that $(u^i, v^i)$ satisfy

\be\label{eq3-20}
\begin{cases}-\Delta u + a(P^i) u =\chi^i(
\mu_1 u^3+\beta uv^2), & x\in \RN,\\
-\Delta v +b(P^i) v =\chi^i(\mu_2 v^3+\beta vu^2),    & x\in
\RN,\\
u(x), v(x)\in H^1(\R^3).\end{cases}\ee

Define
{\allowdisplaybreaks
\begin{align}\label{eq3-21}
\widetilde{L}_i(u, v):=&\frac{1}{2}\intRN{(|\nabla u|^2+a(P^i) u^2+|\nabla v|^2+b(P^i) v^2)}\nonumber\\
&-\frac{1}{4}\int_{\R^3}\chi^i(\mu_1 u^4+2\beta u^2v^2 +\mu_2 v^4)\,dx,
\end{align}
}%
then we see from (\ref{eq2-3}) that
$$\widetilde{L}_i(u^i, v^i)=\max_{t>0}\widetilde{L}_i(t u^i, t v^i)\ge \max_{t>0} L_{P^i}(t u^i, t v^i)\ge m(P^i).$$
Define \be\label{eq30}H_n:=\frac{|\nabla u_{\e_n}|^2+a_{\e_n}|u_{\e_n}|^2+|\nabla v_{\e_n}|^2+b_{\e_n}|v_{\e_n}|^2}{2}
-G_{\e_n}(x, u_{\e_n}, v_{\e_n}).\ee
Then,
{\allowdisplaybreaks
\begin{align}\label{eq3-22}
\liminf_{R\to\iy}\liminf_{n\to\iy} &\int_{B(P_n^i, R)}H_n\,dx\nonumber\\
=\liminf_{R\to\iy}\liminf_{n\to\iy} &\int_{B(0, R)}\Big(
\frac{|\nabla u_n^i|^2+a(\e_n x+\e_n P_n^i)|u_n^i|^2}{2}\nonumber\\
 +&\frac{|\nabla v_n^i|^2+b(\e_n x+\e_n P_n^i)|v_n^i|^2}{2}-G_{\e_n}(x+P_n^i, u_{n}^i, v_{n}^i)\Big)\,dx\nonumber\\
=\widetilde{L}_i(u^i, v^i)&\ge m(P^i).
\end{align}
}%
Similarly, we have
{\allowdisplaybreaks
\begin{gather}\label{eq3-23}
    \limsup_{R\to\iy}\limsup_{n\to\iy}\int_{B(P_n^i, R)\backslash B(P_n^i, R/2)}(|\nabla u_{\e_n}|^2+a_{\e_n}|u_{\e_n}|^2)\,dx=0,\\
\label{eq3-24}\limsup_{R\to\iy}\limsup_{n\to\iy}\int_{B(P_n^i, R)\backslash B(P_n^i, R/2)}(|\nabla v_{\e_n}|^2+b_{\e_n}|v_{\e_n}|^2)\,dx=0.
\end{gather}
}%
Define $B_{R, n}:=\R^3\backslash\cup_{i=1}^k B(P_n^i, R)$. We claim that
\be\label{eq3-25}\liminf_{R\to\iy}\liminf_{n\to\iy} \int_{B_{R, n}}H_n\,dx\ge 0.\ee
From (\ref{eq3-8}) and (\ref{eq3-9}), one has
{\allowdisplaybreaks
\begin{align}\label{eq3-26}\int_{B_{R, n}}H_n\,dx\ge &\frac{1}{2}\int_{B_{R,n}}
\Big(|\nabla u_{\e_n}|^2+a_{\e_n}|u_{\e_n}|^2+|\nabla v_{\e_n}|^2+b_{\e_n}|v_{\e_n}|^2\nonumber\\
&-\nabla_{(u, v)}G_{\e_n}(x, u_{\e_n}, v_{\e_n})(u_{\e_n}, v_{\e_n})\Big)\,dx:= A_{R, n}. \end{align}
}%
Let $\vp_{R, n}\in C_0^{\iy}(\R^3)$ satisfy $\vp_{R,n}=1$ on $B_{R, n}$, $\vp_{R,n}=0$ on $\cup_{i=1}^k B(P_n^i, R/2)$
and $|\nabla \vp_{R, n}|\le C/R$.
Recall that $\nabla_{(u, v)}G_{\e_n}(x, u_{\e_n}, v_{\e_n})(u_{\e_n}, v_{\e_n})\ge 0$.
Multiplying (\ref{eq3-12}) with $(\vp_{R, n}u_{\e_n}, \vp_{R,n}v_{\e_n})$ and integrating over $\R^3$, we have
{\allowdisplaybreaks
\begin{align*}
2A_{R,n}\ge &-\int_{\cup_{i=1}^k B(P_n^i, R)}(|\nabla u_{\e_n}|^2+a_{\e_n}|u_{\e_n}|^2+|\nabla v_{\e_n}|^2+b_{\e_n}|v_{\e_n}|^2)\vp_{R,n}\,dx\\
 &-\int_{\cup_{i=1}^k B(P_n^i, R)}(u_{\e_n}\nabla u_{\e_n}\nabla \vp_{R, n}+v_{\e_n}\nabla v_{\e_n}\nabla \vp_{R, n}):=A^1_{R, n}+A^2_{R, n}.
\end{align*}
}%
From $(\ref{eq3-23})$ and $(\ref{eq3-24})$ we see that $\liminf\limits_{R\to\iy}\liminf\limits_{n\to\iy}A^1_{R, n}=0.$
Since $\cup_{i=1}^k B(P_n^i, R)\subset O_{\e_n}^{2\dd}$ for $n$ large enough, we see from (\ref{eq3-2}) and Lemma \ref{lemma5} that
{\allowdisplaybreaks
\begin{align*}
&\limsup_{R\to\iy}\limsup_{n\to\iy}\int_{\cup_{i=1}^k B(P_n^i, R)}|u_{\e_n}\nabla u_{\e_n}\nabla \vp_{R, n}|\,dx\\
\le &\limsup_{R\to\iy}\limsup_{n\to\iy}C/R\int_{\cup_{i=1}^k B(P_n^i, R)}|\nabla u_{\e_n}|^2+a_\e|u_{\e_n}|^2\,dx\\
\le & \limsup_{R\to\iy}C/R=0.
\end{align*}
}%
Therefore, $\liminf\limits_{R\to\iy}\liminf\limits_{n\to\iy}A^2_{R, n}=0.$
That is, (\ref{eq3-25}) holds. By (\ref{eq3-22}) and (\ref{eq3-25}) we have
$$\liminf\limits_{n\to\iy} J_{\e_n}(u_{\e_n}, v_{\e_n})\ge\sum_{i=1}^k \widetilde{L}_i(u^i, v^i)\ge \sum_{i=1}^k m(P^i).$$
This completes the proof.\hfill$\square$

\bl\label{lemma8}Let $(u_\e, v_\e)$ be in Lemma \ref{lemma3}. Let $P_\e\in O_\e$ such that
$\liminf\limits_{\e\to0}(u_{\e}+v_{\e})(P_\e)>0$. Then $$\liminf\limits_{\e\to 0}\hbox{dist}(\e P_\e, \partial O)>0,\quad\liminf_{\e\to 0\atop{R\to\iy}}\|u_\e+v_\e\|_{L^{\iy}(O_\e\backslash B(P_\e, R))}=0.$$\el

\noindent{\bf Proof. } Assume that there exists $\e_n\to 0$ such that
$\lim\limits_{n\to \iy}\hbox{dist}(\e_n P_{\e_n}, \partial O)=0$. Passing to a subsequence, we may assume that
$\e_n P_{\e_n}\to P_0\in\partial O$. By Lemmas \ref{lemma5} and \ref{lemma7} we have
$$\inf_{p\in O}m(P)\ge \lim_{n\to\iy}J_{\e_n}(u_{\e_n}, v_{\e_n})\ge m(P_0)\ge\inf_{P\in\partial O}m(P),$$
a contradiction with assumption $(V_4)$. Therefore, $\liminf\limits_{\e\to 0}\hbox{dist}(\e P_\e, \partial O)>0$.

Assume that there exists $y_n\in O_{\e_n}$ such that
$$|y_n-P_{\e_n}|\to +\iy,\quad \liminf_{n\to\iy}(u_{\e_n}+v_{\e_n})(y_n)>0.$$
Passing to a subsequence, we may assume that $\e_n y_n\to y_0\in \overline{O}$ and $\e_n P_{\e_n}\to P_0\in O$.
Then by Lemma \ref{lemma7} again, we obtain
$$\inf_{p\in O}m(P)\ge \lim_{n\to\iy}J_{\e_n}(u_{\e_n}, v_{\e_n})\ge m(P_0)+m(y_0)\ge 2\inf_{P\in O}m(P),$$
a contradiction. This completes the proof.\hfill$\square$\\

By Lemmas \ref{lemma6} and \ref{lemma8}, there exists $x_\e\in O_\e$ such that
\be\label{eq3-27}(u_\e+v_\e)(x_\e)=\max_{x\in O_\e}(u_\e+v_\e)(x)\ge\sqrt{\min\{a_0, b_0\}/\bb}.\ee
The following lemma plays a crucial role in the proof of decay estimates.

\bl\label{lemma9}Let $(u_\e, v_\e)$ be in Lemma \ref{lemma3} and $x_\e$ in (\ref{eq3-27}). Let
$(\e_n)_{n\ge 1}$ be any a subsequence with $\e_n\to 0$. Then passing to a subsequence,
$\e_n x_{\e_n}\to P_0\in\m$ and $(u_{\e_n}(x+x_{\e_n}), v_{\e_n}(x+x_{\e_n}))$ converges to some $(U, V)\in \mathcal{S}(P_0)$ strongly
in $C^1_{loc}(\R^3)$, where $S(P_0)$ is defined in (\ref{eq2-11}). Moreover,
\be\label{eq3-28}\lim_{n\to\iy}\int_{\R^3}|\nabla u_{\e_n}(\cdot+x_{\e_n})-\nabla U|^2=0,\,\,\,\lim_{n\to\iy}\int_{\R^3}|\nabla v_{\e_n}(\cdot+x_{\e_n})-\nabla V|^2=0.\ee
In particular, both $u_\e>0$ and $v_\e>0$ hold for $\e>0$ small enough.\el

\noindent{\bf Proof. } By the proof of Lemma \ref{lemma8}, it it easy to see that $\e_n x_{\e_n}\to P_0\in\m$. Repeating the proof of Lemma
\ref{lemma7}, one has that $(u_n(x), v_n(x)):=(u_{\e_n}(x+x_{\e_n}), v_{\e_n}(x+x_{\e_n}))$ converges to some $(U, V)\in H$ strongly
in $C^1_{loc}(\R^3)$. Since $P_0\in\m\in O$, we have that $\chi_{O_{\e_n}}(\cdot+x_{\e_n})$ converges almost everywhere
to $1$. Therefore, $(U, V)$ is a nontrivial solution of (\ref{eq4}) with $P=P_0$. By Lemmas \ref{lemma5} and \ref{lemma7},
$$m(P_0)\ge \limsup\limits_{n\to\iy} J_{\e_n}(u_{\e_n}, v_{\e_n})\ge\liminf\limits_{n\to\iy} J_{\e_n}(u_{\e_n}, v_{\e_n})\ge L_{p_0}(U, V)\ge m(P_0).$$
Therefore, $\lim\limits_{n\to\iy} J_{\e_n}(u_{\e_n}, v_{\e_n})=L_{p_0}(U, V)= m(P_0)$. By Corollary \ref{corollary1}, $U, V>0$.
By (\ref{eq3-27}), $0\in\R^3$ is a maximum point of $U+V$. Combining this with Remark \ref{remark},
we see that $U, V$ are radially symmetric, that is, $(U, V)\in \mathcal{S}(P_0)$.

We claim that
\be\label{eq3-29}\lim_{R\to\iy}\lim_{n\to\iy}\int_{O_{\e_n}\backslash B(x_{\e_n}, R)}
\nabla_{(u, v)}G_{\e_n}(x, u_{\e_n}, v_{\e_n})(u_{\e_n}, v_{\e_n})\,dx=0.\ee
From (\ref{eq3-8}), (\ref{eq3-9}) and Fatou Lemma we have
{\allowdisplaybreaks
\begin{align*}
L_{P_0}(U, V)&=\lim\limits_{n\to\iy} J_{\e_n}(u_{\e_n}, v_{\e_n})\\
=&
\lim_{n\to\iy}\intRN{\frac{1}{2}\nabla_{(u, v)}G_{\e_n}(x, u_{\e_n}, v_{\e_n})(u_{\e_n}, v_{\e_n})-G_{\e_n}(x, u_{\e_n}, v_{\e_n})}\\
\ge&\lim_{n\to\iy}\int_{O_{\e_n}}{\frac{1}{2}\nabla_{(u, v)}G_{\e_n}(x, u_{\e_n}, v_{\e_n})(u_{\e_n}, v_{\e_n})-G_{\e_n}(x, u_{\e_n}, v_{\e_n})}\,dx\\
=&\lim_{n\to\iy}\int_{O_{\e_n}\backslash B(x_{\e_n}, R)}F(u_{\e_n}, v_{\e_n})\,dx+\lim_{n\to\iy}\int_{B(x_{\e_n}, R)}F(u_{\e_n}, v_{\e_n})\,dx\\
=&\lim_{n\to\iy}\int_{(O_{\e_n}-x_{\e_n})\backslash B(0, R)}F(u_{n}, v_{n})\,dx+\lim_{n\to\iy}\int_{B(0, R)}F(u_{n}, v_{n})\,dx\\
\ge &\int_{\R^3\backslash B(0, R)}F(U, V)\,dx+\int_{B(0, R)}F(U, V)\,dx=L_{P_0}(U, V).
\end{align*}
}%
This implies that $\lim\limits_{n\to\iy}\int_{O_{\e_n}\backslash B(x_{\e_n}, R)}F(u_{\e_n}, v_{\e_n})\,dx=\int_{\R^3\backslash B(0, R)}F(U, V)\,dx$.
By (\ref{eq3-8}) again, we see that (\ref{eq3-29}) holds.

Next, we claim that
\be\label{eq3-30}\lim_{R\to\iy}\lim_{n\to\iy}\int_{\R^3\backslash B(x_{\e_n}, R)}|\nabla u_{\e_n}|^2\,dx=0.\ee
Assume by contradiction that, up to a subsequence,
\be\label{eq3-31}\lim_{R\to\iy}\lim_{n\to\iy}\int_{\R^3\backslash B(x_{\e_n}, R)}|\nabla u_{\e_n}|^2\,dx=2\al_2>0.\ee
Let $H_n$ be in (\ref{eq30}). Since $\lim\limits_{n\to\iy} J_{\e_n}(u_{\e_n}, v_{\e_n})=L_{p_0}(U, V)$,
by repeating the proof of Lemma \ref{lemma7} (especially see (\ref{eq3-22}) and (\ref{eq3-25})),
we deduce that
\be\label{eq3-32}\lim_{R\to\iy}\lim_{n\to\iy} \int_{\R^3\backslash B(x_{\e_n}, R)}H_n\,dx= 0.\ee
On the other hand, by (\ref{eq3-9}) and Lemma \ref{lemma5} we have
{\allowdisplaybreaks
\begin{align*}
&\lim_{n\to\iy}\int_{\R^3\backslash O_{\e_n}}
\nabla_{(u, v)}G_{\e_n}(x, u_{\e_n}, v_{\e_n})(u_{\e_n}, v_{\e_n})\,dx\\
\le &\lim_{n\to\iy}\int_{\R^3\backslash O_{\e_n}}\sqrt{\bb}\g_{\e_n}(|x|)(u_{\e_n}^2+v_{\e_n}^2)\,dx\\
\le&\lim_{n\to\iy}\frac{\sqrt{\bb}\e_n^2}{\log{\rho_0/\e_n}}\int_{\R^3\backslash O_{\e_n}}\frac{u_{\e_n}^2+v_{\e_n}^2}{|x|^2}\,dx\\
\le&\lim_{n\to\iy}\frac{4\sqrt{\bb}\e_n^2}{\log{\rho_0/\e_n}}\int_{\R^3}|\nabla u_{\e_n}|^2+|\nabla v_{\e_n}|^2\,dx=0.
\end{align*}
}%
Combining this with (\ref{eq3-29}), we get
\be\label{eq3-33}\lim_{R\to\iy}\lim_{n\to\iy}\int_{\R^3\backslash B(x_{\e_n}, R)}
\nabla_{(u, v)}G_{\e_n}(x, u_{\e_n}, v_{\e_n})(u_{\e_n}, v_{\e_n})\,dx=0.\ee
By (\ref{eq3-26}), (\ref{eq3-31}) and (\ref{eq3-33}) we deduce that
$$\liminf_{R\to\iy}\liminf_{n\to\iy} \int_{\R^3\backslash B(x_{\e_n}, R)}H_n\,dx\ge \al_2>0,$$
which contradicts with (\ref{eq3-32}). Therefore, (\ref{eq3-30}) holds, that is,
$$\lim_{R\to\iy}\lim_{n\to\iy}\int_{\R^3\backslash B(0, R)}|\nabla u_{n}|^2\,dx=0.$$
Since $u_n\to U$ strongly in $C^1_{loc}(\R^3)$, we have $\lim\limits_{n\to\iy}\int_{B(0, R)}|\nabla u_{n}-\nabla U|^2\,dx=0$ for any
$R>0$. Therefore, $\lim\limits_{n\to\iy}\int_{\R^3}|\nabla u_{n}-\nabla U|^2\,dx=0$. Similarly, $\lim\limits_{n\to\iy}\int_{\R^3}|\nabla v_{n}-\nabla V|^2\,dx=0$, and so (\ref{eq3-28}) holds. This means that both $u_\e\not\equiv 0$ and $v_\e\not\equiv0$ for $\e>0$ small enough. By the strong
maximum principle, we have $u_\e, v_\e>0$. This
completes the proof.\hfill$\square$\\

In order to prove that $(u_\e, v_\e)$ is actually a solution of the original problem (\ref{eq3}), we need to give decay estimates of $(u_\e, v_\e)$. First let us recall the following classical result of elliptic estimates.

\bl\label{lemma8-17}(see \cite[Lemma 8.17]{GT}) Let $\Omega$ is an open subset
of $\RN$ and $c\in L^{\iy}(\Omega)$. Suppose that $t>N$, $h\in L^{\frac{t}{2}}(\Omega)$ and
$u\in H^1(\Omega)$ satisfies $-\Delta u(y)+c(y)u(y)\le h(y),\,y\in\Omega$ in
the weak sense. Then for any ball $B(y, 2r)\subset\Omega$,
$$\sup_{B(y, r)}u\le C(\|u^+\|_{L^6 (B(y, 2r))}+\|h\|_{L^{t/2}(B(y,2r))}),$$
where $C= C(N, t, r, \|c\|_{L^{\iy}(\Omega)})$ is independent of $u$ and $y$, and $u^+=\max\{0,
u\}$.\el

\bl\label{lemma10} Let $(u_\e, v_\e)$ be in Lemma \ref{lemma3} and $x_\e$ in (\ref{eq3-27}). Then for $\e>0$ sufficiently small,
there exist some $c, C>0$
independent of $\e>0$, such that
\be\label{eq3-37}\om_\e(x):=u_\e(x)+v_\e(x)\le C \exp\left(-c \,\hbox{dist}(x, \partial O_\e^{3\dd}\cup\{x_\e\})\right),\quad x\in O_\e^{3\dd}.\ee\el

\noindent{\bf Proof. } By (\ref{eq3-5}) and (\ref{eq3-12}) we have
$$\begin{cases}-\Delta u_\e +a_\e u_\e \le \mu_1 u_\e^3+\bb u_\e v_\e^2, & x\in \R^3,\\
-\Delta v_\e +b_\e v_\e \le\mu_2 v_\e^3+\bb u_\e^2 v_\e,    & x\in
\R^3.\end{cases}$$
Without loss of generality, we assume that $a_0\le b_0$. Then by (\ref{eq3-2}) we get
\be\label{eq3-34}-\Delta \om_\e+\frac{a_0}{2}\om_\e\le\bb \om_\e^3 \quad\hbox{in $O_\e^{5\delta}$.}\ee
Note that $6=2^\ast$ and $2>N/2$ in dimension $N=3$.
Then by Lemma \ref{lemma8-17}, there exists $C>0$ independent of small $\e>0$ such that
\be\label{eq3-35}\sup_{x\in B(y, 1)}\om_\e(x)\le  C\left(\|\om_\e\|_{L^6(B(y, 2))}+\|\om_\e\|^3_{L^6(B(y, 2))}\right),\quad\forall\,\, y\in O_\e^{4\dd}.\ee
Since $u_\e, v_\e$ are uniformly bounded in $L^6(\RN)$, we see that $\{\|\om_\e\|_{L^{\iy}(O_\e^{4\dd})}\}_\e$ is
uniformly bounded. Besides,
by (\ref{eq2-1-3}) in Lemma \ref{lemma1}, for any $\sigma>0$, there exists $R>0$ large enough, such that
$$\|U\|_{L^{6}(\R^3\backslash B(0, R))}\le \sg,\quad \|V\|_{L^{6}(\R^3\backslash B(0, R))}\le \sg,\quad\forall\,\, (P,U,V)\in\mathcal{S}.$$
By (\ref{eq3-28}) in Lemma \ref{lemma9} and Sobolev inequalities we deduce that
$$\|u_\e\|_{L^{6}(\R^3\backslash B(x_\e, R))}\le 2\sg,\quad \|v_\e\|_{L^{6}(\R^3\backslash B(x_\e, R))}\le 2\sg,\quad\hbox{for $\e>0$ small enough}.$$
Combining this with (\ref{eq3-35}), one has that
$$\sup_{y\in O_\e^{4\dd}\backslash B(x_\e, R+2)}\om_\e(y)\le C\sg,\quad\hbox{for $\e>0$ small enough}.$$
Then by (\ref{eq3-34}), there is a small $\sg>0$ and so a large $R>0$, such that
\begin{align}\label{eqeq}-\Delta \om_\e+\frac{a_0}{4}\om_\e\le 0\,\,\, \hbox{in $O_\e^{4\dd}\backslash B(x_\e, R)$ and}\,\,\,\sup_{y\in\overline{O_\e^{4\dd}\backslash B(x_\e, R)}}\om_\e(y)\le\mathcal{C}
\end{align}
hold for any $\e>0$ sufficiently small. Here $\mathcal{C}>0$ is independent of $\e$.

Now we want to apply a comparison principle to obtain (\ref{eq3-37}). For any $y\in \partial O^{4\dd}$, we define
open sets $V_y$ as
$$V_y:=\left\{x\in O^{4\dd} \,\,:\,\, |x-y|<\frac{10}{9}\text{dist}(x, \partial O^{4\dd})\right\}.$$
By the finite covering theorem, there exist $m\in\mathbb{N}$ and $y_i\in \partial O^{4\dd}$, $1\le i\le m$, such that
$$\overline{O^{\frac{31}{10}\dd}}\subset\cup_{i=1}^m V_{y_i}.$$
Define a comparison function
$$f_\e(x):=\mathcal{C}e^{\frac{\sqrt{a_0}}{2}R}e^{-\frac{\sqrt{a_0}}{2}|x-x_\e|}
+\mathcal{C}e^{\dd\frac{\sqrt{a_0}}{2\e}}\sum_{i=1}^m e^{-\frac{\sqrt{a_0}}{2}|x-y_i/\e|},$$
then it is easy to check that
$$-\Delta f_\e+\frac{a_0}{4}f_\e>0\quad\text{in}\,\,\,O_\e^{\frac{31}{10}\dd}\setminus B(x_\e, R).$$
Moreover, $f_\e>\mathcal{C}$ on $\partial B(x_\e, R)$. For any $x\in\partial O_\e^{\frac{31}{10}\dd}$, $\e x\in V_{y_j}$ for some $1\le j\le m$ and so
$$|\e x-y_j|<\frac{10}{9}\text{dist}(\e x, O^{4\dd})=\dd,$$
which implies
$f_\e(x)>\mathcal{C}e^{\dd\frac{\sqrt{a_0}}{2\e}}e^{-\frac{\sqrt{a_0}}{2}|x-y_j/\e|}\ge \mathcal{C}$. Combining these with (\ref{eqeq}), we deduce from the maximum principle that
$$w_\e(x)\le f_\e(x),\quad\forall\,x\in O_\e^{\frac{31}{10}\dd}\setminus B(x_\e, R).$$
For any $x\in O_\e^{3\dd}\setminus B(x_\e, R)$, we have $|x-y_i/\e|\ge \text{dist}(x, \partial O_\e^{4\dd})=\text{dist}(x, \partial O_\e^{3\dd})+\dd/\e$ for all $1\le i\le m$, so
{\allowdisplaybreaks
\begin{align*}
w_\e(x)&\le f_\e(x)=\mathcal{C}e^{\frac{\sqrt{a_0}}{2}R}e^{-\frac{\sqrt{a_0}}{2}|x-x_\e|}
+\mathcal{C}e^{\dd\frac{\sqrt{a_0}}{2\e}}\sum_{i=1}^m e^{-\frac{\sqrt{a_0}}{2}|x-y_i/\e|}\\
&\le\mathcal{C}e^{\frac{\sqrt{a_0}}{2}R}e^{-\frac{\sqrt{a_0}}{2}|x-x_\e|}
+m\mathcal{C}e^{-\frac{\sqrt{a_0}}{2}\text{dist}(x, \partial O_\e^{3\dd})}\\
&\le(\mathcal{C}e^{\frac{\sqrt{a_0}}{2}R}
+m\mathcal{C})\exp\left(-\frac{\sqrt{a_0}}{2}\hbox{dist}\left(x, \partial O_\e^{3\dd}\cup\{x_\e\}\right)\right).
\end{align*}
}%
That is, (\ref{eq3-37}) holds.
This completes the proof.
\hfill$\square$\\

By $(V_2)$ there exists $R_1>0$ large enough, such that $O\subset B(0, R_1)$ and for some $c>0$,
\be\label{eq3-38}a(x), \,\,b(x)\ge \frac{c}{|x|^2\log (|x|)}, \quad\forall \,\,|x|\ge R_1.\ee

\bl\label{lemma11} Let $(u_\e, v_\e)$ be in Lemma \ref{lemma3} and $x_\e$ in (\ref{eq3-27}). Then for sufficiently large $R_2>R_1$,
there exists $c, C>0$
independent of $\e>0$, such that
\be\label{eq3-39}\om_\e(x)\le C e^{-\frac{c}{\e}}, \quad\hbox{for all $\dd/\e\le|x-x_\e|\le 2R_2/\e$}\ee
holds for $\e>0$ sufficiently small.\el

\noindent{\bf Proof. } Let $D:=\{x\in\R^3 : a(x)=0\,\,\hbox{or}\,\, b(x)=0\}$. Then by $(V_2)-(V_3)$ and (\ref{eq3-38}) we see that $D\subset\R^3\backslash O$ is compact and $D\subset B(0, R_1)$.
First, we claim that, for any sufficiently large $R_2>R_1$ and sufficiently small $l>0$, there exists $C, c>0$ independent of small $\e>0$, such that
\be\label{eq3-40}\om_\e(x)\le C e^{-\frac{c}{\e}},\,\,\hbox{for $\dd/\e\le|x-x_\e|\le 2R_2/\e$ and dist$(\e x, D)\ge l$}.\ee
By the definition of $D$, we may assume that
$$\inf_{x\in B(x_\e, 5R_2/\e)\backslash D^{l/4}_\e}\min\{a_\e(x), b_\e(x)\}=a'>0.$$
By a similar proof of Lemma \ref{lemma10}, we can prove that
$$-\Delta \om_\e+\frac{a'}{4}\om_\e\le 0\quad\text{and}\,\,\,\om_\e(x)\le\mathcal{C}'$$
hold for $\dd/\e\le|x-x_\e|\le 4R_2/\e$ and dist$(\e x, D)\ge l/3$ when $\e>0$ small enough.
Recall from (\ref{eq3-37}) that $w_\e(x)\le Ce^{-\frac{c}{\e}}$ for any $|x-x_\e|=\dd/\e$.
By a similar proof of Lemma \ref{lemma10} (i.e., consider $B(x_\e, 4R_2/\e)\setminus D^{l/3}_\e$ and $B(x_\e, 3R_2/\e)\setminus D^{l/2}_\e$ similar as $O^{4\dd}_\e$ and $O^{3\dd}_\e$ in Lemma \ref{lemma10} respectively), there exist some $C, c>0$ independent of small $\e>0$, such that
\be\label{eq3-41}\om_\e(x)\le C \exp\left(-c \,\hbox{dist}\left(x, \partial D_\e^{l/2}\cup\partial B(x_\e, 3R_2/\e)\cup\{x_\e\}\right)\right)\ee
holds for $\dd/\e\le|x-x_\e|\le 3R_2/\e$ and dist$(\e x, D)\ge l/2$. That is, (\ref{eq3-40}) holds.

Let $l>0$ small enough such that $D^{2l}\cap O=\emptyset$. Let $\psi\ge 0$ satisfy
$$\begin{cases}-\Delta \psi =\la_1 \psi, & x\in D^{2l},\\
\psi=0,    & x\in
\partial D^{2l},\end{cases}$$
where $\la_1$ is the first eigenvalue. We may assume that $\max\limits_{x\in D^{2l}}\psi(x)=1$. Define $\psi_\e(x):=\psi(\e x)$.
By (\ref{eq3-15}) we see that
{\allowdisplaybreaks
\begin{align*}
&-\Delta \psi_\e +a_\e \psi_\e -\frac{\partial_u G_\e(x, u_\e, v_\e)}{u_\e}\psi_\e\ge\la_1 \e^2\psi_\e-\al_0\g_\e(|x|)\psi_\e\\
\ge &\e^2\left(\la_1-\frac{\al_0}{|\rho_0/\e|^2\log{|\rho_0/\e|}}\right)\psi_\e\ge 0 \,\,\,\hbox{in $D_\e^{l}$, for $\e>0$ small enough.}
\end{align*}
}%
Recall from (\ref{eq3-40}) that $u_\e(x)\le w_\e(x)\le Ce^{-\frac{c}{\e}}$ for all $x\in\partial D_\e^l$. Besides, there exists $c>0$ such that $\min\limits_{x\in \partial D^{l}}\psi(x)=c$.
Again by a comparison principle, there exist some $C, c>0$ independent of small $\e>0$, such that
$$u_\e(x)\le C e^{-\frac{c}{\e}} \psi_\e (x)\le Ce^{-\frac{c}{\e}}\,\,\,\hbox{for $x\in D_\e^l$.}$$
Similarly, we can prove that
$v_\e(x)\le C e^{-\frac{c}{\e}}\,\,\,\hbox{for $x\in D_\e^l$.}$ Therefore,
$\om_\e(x)\le C e^{-\frac{c}{\e}}\,\,\,\hbox{for $x\in D_\e^l$.}$ Combining this with (\ref{eq3-40}), we see that (\ref{eq3-39}) holds.
This completes the proof.\hfill$\square$

\bl\label{lemma12}Let $(u_\e, v_\e)$ be in Lemma \ref{lemma3}, $x_\e$ in (\ref{eq3-27}) and $R_2$ in Lemma \ref{lemma11}.
Then for any $\al>0$, there exists $c, C>0$
independent of $\e>0$, such that
\be\label{eq3-42}\om_\e(x)\le C e^{-\frac{c}{\e}}|x|^{-1}|\log{|x|}|^{-\al}, \quad\hbox{for all $x\in\R^3\backslash B(0,R_2/\e)$}\ee
holds for $\e>0$ sufficiently small. In particular, there exists $\e_0>0$ small enough, such that for any $\e\in (0,\e_0)$,
$(u_\e, v_\e)$ is a positive vector solution of (\ref{eq3-1}).\el

\noindent{\bf Proof. } The following proof is similar to \cite{BB}, and we give the proof for the completeness.
For any fixed $\al>0$, we define
\be\label{eq3-43}\G_\e(x)=\frac{1}{|x|(\log|x|)^\al},\ee
then there exists some $C>0$ such that $\min\limits_{x\in \partial O_\e} \G_\e(x)\ge C\e^2$. For any $x\in \R^3\backslash B(0, R_2/\e)$, we have
$$\Delta \G_\e(x)=\frac{\al}{|x|^3(\log|x|)^{\al+1}}+\frac{\al(\al+1)}{|x|^3(\log|x|)^{\al+2}},$$
and so it follows from (\ref{eq3-15}) and (\ref{eq3-38}) that
{\allowdisplaybreaks
\begin{align*}
&\left(-\Delta \G_\e +a_\e \G_\e -\frac{\partial_u G_\e(x, u_\e, v_\e)}{u_\e}\G_\e\right)\Big/\G_\e\\
\ge&\frac{c}{|\e x|^2\log|\e x|}-\frac{\al}{|x|^2\log|x|}-\frac{\al(\al+1)}{|x|^2(\log|x|)^{2}}-\frac{\al_0 \e^2}{|x|^2\log|x|}.
\end{align*}
}%
Since for small $\e>0$ and $|x|\ge R_2/\e$,
$$\frac{1}{|\e x|^2\log|\e x|}\ge \frac{1}{\e(\e+1)}\frac{1}{| x|^2\log|x|},$$
we see that
$$-\Delta \G_\e +a_\e \G_\e -\frac{\partial_u G_\e(x, u_\e, v_\e)}{u_\e}\G_\e\ge 0\,\,\,\hbox{in $\R^3\backslash B(0, R_2/\e)$}.$$
For $x\in \partial B(0, R_2/\e)$, since $\e x\not\in O$ and $\e x_\e\in \m^\dd$ for $\e>0$ sufficiently small by Lemma \ref{lemma9},
we have
$$2R_2/\e\ge|x-x_\e|\ge\frac{|\e x-\e x_\e|}{\e}\ge\frac{dist(\R^3\backslash O, \m^\dd)}{\e}\ge\frac{4\dd}{\e},$$
by (\ref{eq3-39}) in Lemma \ref{lemma11}, $u_\e(x)\le Ce^{-\frac{c}{\e}}$ for all $x\in \partial B(0, R_2/\e)$. Therefore, by a comparison
principle, there exists $C, c>0$ independent of small $\e>0$ such that
$$u_\e(x)\le C e^{-\frac{c}{\e}}|x|^{-1}|\log{|x|}|^{-\al}, \quad\hbox{for all $x\in\R^3\backslash B(0,R_2/\e)$}$$
holds for $\e>0$ sufficiently small.
By a similar proof, the conclusion also holds for $v_\e$. That is, (\ref{eq3-42}) holds for $\e>0$ sufficiently small.

Now we fix a $\al>1/2$. For any $x\in\R^3\backslash B(0, R_2/\e)$, we have
{\allowdisplaybreaks
\begin{align*}
4F(u_\e(x), v_\e(x))=&\mu_1 u_\e^4(x)+2\bb u_\e^2(x)v_\e^2(x)+\mu_2 v_\e^4(x)\\
\le &\bb \om_\e^4(x)\le C e^{-\frac{4c}{\e}}|x|^{-4}|\log{|x|}|^{-4\al}\\
<& \frac{\e^4}{|x|^4(\log |x|)^2}=\g_\e(|x|)^2\,\,\,\hbox{for $\e>0$ small enough.}
\end{align*}
}%
For $x\in B(0, R_2/\e)\backslash O_\e$, we have $\dd/\e<\rho_0/\e\le|x|\le R_2/\e$. Then by (\ref{eq3-39}) in Lemma \ref{lemma11}, we deduce that
$$4F(u_\e(x), v_\e(x))\le Ce^{-\frac{4c}{\e}}<\frac{\e^4}{|x|^4(\log |x|)^2}=\g_\e(|x|)^2\,\,\hbox{for $\e>0$ small enough.}$$
Therefore, there exists $\e_0>0$ sufficiently small, such that for any $\e\in (0,\e_0)$, we have
$$F(u_\e(x), v_\e(x))<\frac{1}{4}\g_\e(|x|)^2,\,\,\,\forall\,x\in \R^3\backslash O_\e,$$
which implies that $G_\e(x, u_\e, v_\e)\equiv F(u_\e, v_\e)$, and so $(u_\e, v_\e)$ is a positive vector solution of (\ref{eq3-1}). This completes the proof.
\hfill$\square$\\

\noindent {\bf Completion of the proof of Theorem \ref{th1}. } Let $\e\in (0,\e_0)$. Since $x_\e\in O_\e$, we have
$|x_\e|<R_2/\e$ and so
$$|x|>R_2/\e,\quad |x-x_\e|\le 2|x|,\,\,\,\forall\,\,x\in\R^3\backslash B(x_\e, 2R_2/\e).$$
Combining these with (\ref{eq3-42}), there exists $C, c>0$ independent of $\e$ such that
\be\label{eq3-44}\om_\e(x)\le C e^{-\frac{c}{\e}}|x-x_\e|^{-1}|\log{|x-x_\e|}|^{-\al}, \quad\hbox{for all $x\in\R^3\backslash B(x_\e,2R_2/\e)$}.\ee

Define $(\tilde{u}_\e(x), \tilde{v}_\e(x)):=(u_\e(x/\e), v_\e(x/\e))$ and $\tilde{x}_\e:=\e x_\e$. Then $(\tilde{u}_\e, \tilde{v}_\e)$ is
a positive vector solution of (\ref{eq3}). Moreover, $\tilde{x}_\e$ is a maximum point of $\tilde{u}_\e+\tilde{v}_\e$.
Conclusions $(i)$ and $(ii)$ in Theorem \ref{th1} follow directly from Lemma \ref{lemma9}. By (\ref{eq3-44}) we have that
{\allowdisplaybreaks
\begin{align}\label{eq3-45}
\tilde{u}_\e(x)+\tilde{v}_\e(x)&\le C e^{-\frac{c}{\e}}|x/\e-x_\e|^{-1}|\log{|x/\e-x_\e|}|^{-\al}\nonumber\\
&=C e^{-\frac{c}{\e}}\frac{\e}{|x-\tilde{x}_\e||\log(|x-\tilde{x}_\e|/\e)|^{\al}}\nonumber\\
&\le C e^{-\frac{c}{\e}}\frac{1}{|x-\tilde{x}_\e||\log(|x-\tilde{x}_\e|+2)|^{\al}}
\end{align}
}%
holds for all $x\in \R^3\backslash B(\tilde{x}_\e, 2R_2)$. By (\ref{eq3-39}) in Lemma \ref{lemma11}, there exists some $C, c>0$
independent of $\e\in (0,\e_0)$, such that (\ref{eq3-45}) holds for
all $x\in \R^3\backslash B(\tilde{x}_\e, \dd)$.

For any $x\in B(\tilde{x}_\e, \dd)$, since $\tilde{x}_\e=\e x_\e\in\m^\dd$ for $\e>0$ small,
we have $\hbox{dist}(x, \partial O^{3\dd}\cup\{\tilde{x}_\e\})=|x-\tilde{x}_\e|$. By (\ref{eq3-37}) in Lemma \ref{lemma10}, we get that
$$\tilde{u}_\e(x)+\tilde{v}_\e(x)\le C \exp\left({-\frac{c}{\e}|x-\tilde{x}_\e|}\right),\quad x\in  B(\tilde{x}_\e, \dd).$$
Therefore, $(iii)$ in Theorem \ref{th1} holds. This completes the proof.\hfill$\square$

\vskip0.1in

\noindent{\it Acknowledgements.} The authors wish to thank the referee very much for his/her careful reading and valuable comments, which helped to improve some arguments of the paper a lot.

\end{document}